\nonstopmode
\documentclass{amsart}
\usepackage{amssymb}
\usepackage[all]{xy}
\usepackage{verbatim}
\newenvironment{proclaim}[1]{\par\medskip\noindent{\bf #1.}\it}{\par\smallskip}
\newenvironment{demo}[1]{\par\smallskip\noindent{\bf #1.}}{\par\smallskip}

\def\nmb#1#2{#2}         
\def\cit#1#2{\ifx#1!\cite{#2}\else\bibitem{#2}\fi} 
\def\idx{}               
\def\ign#1{}             
\def\o{\circ}

\def\ze{\zeta}

\def\si{\sigma}
\def\ta{\tau}
\def\ph{\varphi}

\def\ps{\psi}
\def\om{\omega}
\def\Ga{\Gamma}

\def\i{^{-1}}

\def\p{\partial}
\def\on{\operatorname}
\begin{document}

\title[Tensor fields on holomorphic orbit spaces of 
    finite groups]{
        Tensor fields and connections on holomorphic orbit spaces of 
     finite groups 
}
\author[Kriegl, 
        Losik, Michor]{ 
        Andreas Kriegl, 
        Mark Losik, and Peter W. Michor  
}
\address{
        A.\  Kriegl: Institut f\"ur Mathematik, Universit\"at Wien,
        Strudlhofgasse 4, A-1090 Wien, Austria
}
\email{ Andreas.Kriegl@univie.ac.at }
\address{
        P.\  W.\  Michor: Institut f\"ur Mathematik, Universit\"at Wien,
        Strudlhofgasse 4, A-1090 Wien, Austria; {\it and:} 
        Erwin Schr\"odinger Institut f\"ur Mathematische Physik,
        Boltzmanngasse 9, A-1090 Wien, Austria
}
\email{ Peter.Michor@esi.ac.at }
\address{
        M. Losik: Saratov State University, ul. Astrakhanskaya, 83,
        410026 Saratov, Russia
}
\email{ LosikMV@info.sgu.ru }
\begin{abstract}
For a representation of a finite group $G$ on a complex vector space 
$V$ we determine when a holomorphic $\binom{p}{q}$-tensor field on 
the principal stratum of the orbit space $V/G$ can be lifted to a 
holomorphic $G$-invariant tensor field on $V$. This extends also to 
connections. As a consequence we determine those holomorphic 
diffeomorphisms on $V/G$ which can be lifted to orbit preserving 
holomorphic diffeomorphisms on $V$. This in turn
is applied to characterize complex orbifolds.
\end{abstract}
\thanks{
        M.L. and P.W.M. were supported  
        by `Fonds zur F\"orderung der wissenschaftlichen Forschung, 
     Projekt P~14195~MAT'
}

\maketitle

\section*{\nmb0{1}. Introduction}
Locally, an orbifold $Z$  can be identified with the orbit space $B/G$,
where $B$ is a $G$-invariant neighborhood of the origin in
a vector space $V$ with a finite group $G\subset GL(V)$ and, 
using this identification,
one can easily define local (and then global)
tensor fields and other differential geometrical objects in $Z$ as 
appropriate  $G$-invariant tensor fields and objects on $B
\subset V$. In particular, one can naturally define Riemannian
orbifolds, Einstein orbifolds, symplectic orbifolds, 
K\"ahler-Einstein orbifolds etc.
 
We study  complex orbifolds, that is, orbifolds 
modeled on orbit spaces $V/G$, where $G$ is 
a finite subgroup of $GL(V)$ for 
a complex vector space $V$. In particular, the orbit 
spaces $Z = M/G$ of a discrete proper group $G$ of 
holomorphic transformations of a complex manifold $M$ are complex 
orbifolds.

An orbifold $X$ has a structure defined by the sheaf 
$\mathfrak F_{X}$ of local invariant holomorphic functions in a local 
uniformizing system. $X$ has also a stratification by strata $S$ 
which are glued from local isotropy type strata of local uniformizing 
systems. In particular, the regular stratum $X_0$ is an open dense 
complex manifold in $X$.
 
Holomorphic geometric objects on $X$ (e.g. tensor fields and 
connections) are  locally defined as invariant objects on the 
uniformizing system. Their restrictions to the regular stratum $X_0$ 
are usual holomorphic  geometric objects on the complex manifold 
$X_0$. 
 
A natural question is to characterize these restrictions, i.e.\ to 
describe tensor fields and connections on $X_0$ which are extendible 
to $X$. We look at the lifting problem for connections because this 
allows a very elegant approach to the lifting problem for holomorphic 
diffeomorphisms. And the last problem has immediate consequences for 
characterizing complex orbifolds, i.e., for answering the following 
question:
Which data does one need besides 
$\mathfrak F_X$ and $X_0$ to characterize a complex orbifold $X$?
The main goal of the paper is to answer these questions. 

We have first to investigate the local situation, thus we consider a 
finite subgroup $G\subset GL(V)$ and the orbit space $Z=V/G$ with the 
structure given by the sheaf $\mathfrak F_{V/G}$ of invariant 
holomorphic functions on $V$, and the orbit type stratification. The 
prime role is played by strata of codimension 1 with the orders of 
the corresponding stabilizer groups, which are arranged in the 
{\it reflection divisor} $D_{V/G}$ which keeps track of all complex 
reflections in $G$.
It turns out that the union $Z_1$ of $Z_0$ and of all codimension 1 
strata is a complex manifold, see \nmb!{3.5}.
We characterize all $G$-invariant holomorphic tensor fields and 
connections on $V$ in terms of the {\it reflection divisor} of the 
corresponding meromorphic tensor field and connection on $Z_1$, see 
\nmb!{3.7} and \nmb!{4.2}. Our result gives a generalization 
\nmb!{3.9} of Solomon's theorem \cite{6}, see \nmb!{3.10}.
Using the lifting property of connections we are able to prove that a 
holomorphic diffeomorphism $Z=V/G\to V/G'=Z'$ between two orbit spaces 
has a holomorphic lift to $V$ which is equivariant over an 
isomorphism $G\to G'$ if and only if $f$ respects the regular strata 
and the reflection divisors, i.e.\ $f(Z_0)\subset Z'_0$ and 
$f_*(D_Z)\subset D_{Z'}$.  In fact we give two proofs of this result, 
which in \cite{LMP} is carried over to the algebraic geometry 
setting for algebraically closed ground fields of characteristic 0. 
The related problem of lifting (smooth) homotopies from (general) 
orbit spaces has been treated in \cite{Bier} and \cite{5}. 

Applying the local results we prove that a complex orbifold $X$ is 
uniquely determined by the sheaf $\mathfrak F_X$, the regular stratum 
$X_0$, and the reflection divisor $D_X$ alone, see \nmb!{6.6}.

\section*{\nmb0{2}. Preliminaries}

\subsection*{\nmb.{2.1}. The orbit type stratification} 
Let $V$ be an $n$-dimensional complex vector space, $G$ a finite subgroup 
of $GL(V)$, and $\pi:V\to V/G$ the quotient projection. The ring  
$\mathbb C[V]^G$ has a minimal system of homogeneous 
generators $\sigma^1,\dots,\sigma^m$. We will use the map 
$\sigma=(\sigma^1,\dots,\sigma^m):V\to\mathbb C^m$. 
Denote by $Z$ the affine algebraic variety in $\mathbb C^m$ defined by 
the relations between $\sigma^1,\dots,\sigma^m$. It is known 
that $\sigma(V)=Z$. 
 
We consider the orbit space $V/G$ endowed with the quotient topology as 
a local ringed space defined by the following sheaf of rings 
$\mathfrak F_{V/G}$: 
if $U$ is an open subset of $V/G$, $\mathfrak F_{V/G}(U)$ is equal 
to the space of $G$-invariant holomorphic 
functions on $\pi^{-1}(U)$. Clearly one may consider 
sections of $\mathfrak F_{V/G}$ on $U$ as functions on $U$. We call these 
functions holomorphic functions on $U$. It is known that the map of the 
orbit space $V/G$ to $Z$ induced by the map $\sigma$ is a homeomorphism.
Moreover, this homeomorphism induces an isomorphism of the sheaf 
$\mathfrak F_{V/G}(U)$ and the structure sheaf of the complex algebraic 
variety $Z$ (see [7]).  Via the above isomorphism we identify the 
local ringed spaces $V/G$ and $Z$. Under this identification the projection 
$\pi$ is identified with the map $\sigma$. 
Let $G$ and $G'$ be finite subgroups of $GL(V)$ and let $Z=V/G$ and $Z'=V/G'$ 
be the corresponding orbit spaces.
By definition a holomorphic diffeomorphism  
of the orbit space $Z$ to the orbit space $Z'$ is an isomorphism of 
$Z$ to $Z'$ as local ringed spaces.  

Let $K$ be a subgroup of $G$, $(K)$ the conjugacy class of $K$. Denote  
by $V_{(K)}$ the set of points of $V$ whose isotropy groups belong to $(K)$  
and put $Z_{(K)}=\pi(V_{(K)})$. It is known that $\{Z_{(K)}\}$ is a 
finite stratification of $Z$, called the isotropy type stratification, 
into locally closed irreducible smooth algebraic subvarieties (see 
\cit!{2}). 
Denote by $Z^i$ the union of the strata of codimension greater than $i$ 
and put $Z_i=Z\setminus Z^i$. Then $Z_0$ is the principal stratum of  
$Z$, i.e.\ $Z_0=Z_{(K)}$ for $K=\{\operatorname{id}\}$. 
It is known that $Z_0$ is a Zariski open subset of $Z$ and a 
complex manifold. It is clear that the restriction 
of the map $\sigma$ to the set $V_{\text{reg}}$ of regular points of $V$ is 
an étale map onto $Z_0$. 
 
In this paper we consider the orbit space $Z=V/G$ with the above structure 
of local ringed space and the stratification $\{Z_{(K)}\}$.   

\subsection*{\nmb.{2.2}. The divisor of a tensor field}
We shall use divisors of meromorphic functions on a complex manifold $X$. 
For technical reasons (see e.g.\ the last formula of this section) 
we define $\on{div}(0)=\sum_S\infty.S$, where the 
sum runs over all complex subspaces of $X$ of codimension $1$.  

Let $f$ and $g$ be two meromorphic functions on $X$.  Then we have 
$\on{div}(f+g)\ge \min\{\on{div}(f),\on{div}(g)\}$, where $\on{div}(f)$ 
denote the divisor of $f$.
Taking the minimum means: For each irreducible complex subspace $S$ of $X$   
of codimension $1$ belonging to the support of $f$ or $g$ 
take the minimum of the coefficients in $\mathbf Z$ of $S$ in $\on{div}(f)$ 
and $\on{div}(g)$. 

Let $P$ be a meromorphic tensor field (i.e., with meromorphic 
coefficient functions in local coordinates) on  
$X$. In local holomorphic coordinates $y^1,\dots,y^n$ on an open 
subset $U\subset X$ the tensor field $P$ can be written as 
$$
P|_U=\sum_{i_1,\dots,i_p,j_1,\dots,j_q} P^{i_1\dots i_p}_{j_1\dots j_q}
\frac{\partial}{\partial y^{i_1}}\otimes\dots\otimes
\frac{\partial}{\partial y^{i_p}}\otimes dy^{j_1}\otimes dy^{j_q}.
$$  
and we define the {\it divisor} of $P$ on $U$ as the minimum of all 
divisors $\on{div}(P^{i_1\dots i_p}_{j_1\dots j_q})\in \on{Div}(U)$ 
for all coefficient functions of $P$. 
The resulting 
coefficient of the complex subspace $S$ of codimension $1$ 
in $\on{div}(P)\in \on{Div}(U)$ 
does not depend on the choice of the holomorphic coordinate system; 
e.g., for a vector field 
$\sum_i X^i\frac{\p}{\p y^i} =  \sum_{i,k} 
X^i\frac{\p u^k}{\p y^i}\frac{\p}{\p u^k}$
we have 
$$
\on{div}\Bigl(\sum_i X^i\frac{\p u^k}{\p y^i}\Bigr) \ge
\min_i\on{div}\Bigl(X^i\frac{\p u^k}{\p y^i}\Bigr) =
\min_i\Bigl(\on{div}(X^i)+\on{div}\Bigl(\frac{\p u^k}{\p y^i}\Bigr)\Bigr) \ge 
\min_i\on{div}(X^i).  
$$

Finally we define the divisor of $P$ on $X$ by gluing the local 
divisors for any holomorphic atlas of $X$. Note that a tensor field 
$P$ is holomorphic if and only if $\on{div}(P)\ge 0$.

\section*{ \nmb0{3}. Invariant tensor fields} 

\subsection*{\nmb.{3.1}} 
Let $P$ be a $G$-invariant holomorphic tensor field of type  
$\binom pq$ on  
$V$. 
Since $\sigma$ is an étale map on $V_{\text{reg}}$, there is a unique  
holomorphic 
tensor field $Q$ on $Z_0$ of type $\binom pq$ such that the pullback  
$\sigma^*(Q)$ 
coincides with the restriction of $P$ to $V_{\text{reg}}$. 
It is clear that the tensor field $P$ is uniquely defined by $Q$. 
 
Consider a holomorphic tensor field $Q$ of type $\binom pq$ on $Z_0$ and its 
pullback $\sigma^*(Q)$ which is a $G$-invariant holomorphic tensor field  
on $V_{\text{reg}}$. Then by the Hartogs extension theorem, $\sigma^*(Q)$  
has a $G$-invariant holomorphic extension to $V$ iff it has a  
holomorphic  
extension to $\sigma^{-1}(Z_1)$. 
 
Denote by $\mathfrak H$ the set of all reflection hyperplanes  
corresponding to all complex reflections in $G$ and, for 
each $H\in\mathfrak H$, by $e_H$ the order of the cyclic subgroup of $G$ 
fixing $H$. It is clear that $\sigma(\cup_{H\in\mathfrak H}H)$ contains all 
strata of codimension $1$. This implies immediately the following 
 
\begin{proclaim}{\nmb.{3.2}. Proposition} 
If $\mathfrak H=\emptyset$, for each  
holomorphic tensor field $P_0$ on $Z_0$ the pullback $\sigma^*(P_0)$  
has a $G$-invariant holomorphic extension to $V$. \qed 
\end{proclaim} 

\subsection*{\nmb.{3.3}. The reflection divisor of the orbit space} 
Consider the set $R_Z$ of all hyper surfaces $\si(H)$ in $Z$, where $H$ runs 
through all reflection hyperplanes in $V$. Note that $\si(H)$ is a 
complex subspace of $Z_1$ of codimension $1$.  
We endow each $S=\si(H)\in R_Z$ with the label $e_H$ of the hyperplane  
$H$. 
It is easily seen that this label does not depend on the choice of  
$H$, we denote it by $e_S$ and we consider $e_S.\,S$ as an effective  
divisor on $Z$ and we consider the effective divisor in $Z_1$ 
$$ 
D=D_{V/G}=D_Z = \sum_{S\in R_Z} e_S.\,S, 
$$ 
which we call the {\it reflection divisor}. 

\subsection*{\nmb.{3.4}. Basic example} 
Let the cyclic group $\mathbb Z_r=\mathbb Z/r\mathbb Z$ with generator  
$\ze_r=e^{2\pi i/r}$ act on $\mathbb C$ by  
$z\mapsto e^{2\pi ik/r}z$ for $r\geq 2$.  
The generating invariant is $\ta(z)=z^r$.   
 
We consider first a holomorphic tensor field  
$P=f(z)(dz)^{\otimes q}\otimes (\tfrac\partial{\partial z})^{\otimes p}$ 
on $\mathbb C$. It is invariant, $\ze_r^*P=P$, if and only if  
$f(\ze_r z)= \ze_r^{p-q}f(z)$, so that in the expansion  
$f(z)=\sum_{k\ge 0}f_kz^k$ at 0 of $f$ the coefficient $f_k\ne 0$  
at most when $k\cong p-q \mod r$.  
Writing $p-q=rs+t$ with $s\in\mathbb Z$ and $0\leq t<r$ 
we see that $P$ is invariant if and only if  
$f(z)=z^{t}g(z^r)$ for holomorphic $g$.  
 
We use the coordinate $y=\ta(z)=z^r$ on $\mathbb C/\mathbb Z_r= \mathbb C$,   
$\ta^* dy = rz^{r-1}dz$ and  
$\ta^*(\tfrac\partial{\partial y}|_{\mathbb C\setminus 0})= 
\frac{1}{r z^{r-1}}\,\frac{\partial}{\partial z}|_{\mathbb C\setminus 0}$,  
and we write  
\begin{align*} 
P|_{\mathbb C\setminus 0}&=g(z^r)z^t(dz)^{\otimes q}\otimes (\tfrac\partial{\partial z})^{\otimes p}\\ 
&=g(y)z^t(rz^{r-1})^{p-q}(dy)^{\otimes q}\otimes (\tfrac\partial{\partial y})^{\otimes p}\\ 
&=g(y)z^{-rs}(rz^r)^{p-q}(dy)^{\otimes q}\otimes (\tfrac\partial{\partial y})^{\otimes p}\\ 
&=g(y)r^{p-q}y^{p-q-s}(dy)^{\otimes q}\otimes (\tfrac\partial{\partial y})^{\otimes p} 
\end{align*} 
(we omitted $\ta^*$).
Thus a holomorphic tensor field $P$ of type $\binom pq$ on  
$\mathbb C$ is $\mathbb Z_r$-invariant if and only if  
$P|_{\mathbb C\setminus 0}=\ta^*Q$ for a meromorphic tensor field
\begin{displaymath} 
Q=g(y)y^m(dy)^{\otimes q}\otimes (\tfrac\partial{\partial y})^{\otimes p} 
\end{displaymath} 
on $\mathbb C$ with  
$g$ holomorphic with $g(0)\ne 0$ and with  
\begin{equation*} 
m\geq p-q-s. 
\end{equation*} 
It is easily checked that the above inequality is equivalent to the following 
one 
\begin{equation*} 
mr+(q-p)(r-1)\ge 0. 
\end{equation*}  
\subsection*{\nmb.{3.5}} 
Suppose $\mathfrak H\neq\emptyset$. Let $z\in Z_1\setminus Z_0$ and 
$v\in\sigma^{-1}(z)$. Then there is a unique hyperplane $H\in\mathfrak H$ such 
that $v\in H$ and the isotropy group $G_v$ is isomorphic to a cyclic  
group. It is evident that the order $r_z=e_H$ of $G_v$ depends only on  
$z=\si(v)$ and is locally constant on $Z_1\setminus Z_0$. 
 
By the holomorphic slice theorem 
(see \cit!{2}, \cit!{3}) there is a  
$G_v$-invariant open 
neighborhood $U_v$ of $v$ in $V$ such that  
the induced map $U_v/G_v\to V/G$ is a local biholomorphic map at $v$. 
 
Choose orthonormal coordinates $z^1,\dots,z^n$ in $V$ with respect to 
a $G$-invariant Hermitian inner product on $V$,  
so that $H=\{z^n=0\}$.  
Then the ring $\mathbb C[V]^{G_v}$ is generated by  
$z^1,\dots,z^{n-1},(z^n)^r$, where 
$r=r_z$. 
 
Put $\tau^1=z^1,\dots,\tau^{n-1}=z^{n-1}$, $\tau^n=(z^n)^r$, 
and $\tau=(\tau^1,\dots,\tau^n):U_v\to\mathbb C^n$.  
Then there are holomorphic functions 
$f^i$ $(i=1,\dots,n)$ in an open neighborhood $W_z$ of $z\in\mathbb C^m$  
such that $\tau^a=f^a\o\sigma|_{U_v}$. 
On the other hand, we know that in an open neighborhood of $v$ all 
$\sigma^a$ for $(a=1,\dots,m)$ are holomorphic functions of the $\tau^i$. 
We denote by $y^i$ the holomorphic function on $Z$ such that 
$\tau^i=y^i\o\sigma$. 
Then we can use $y^i$ as coordinates of $Z$ 
defined in the open neighborhood $W_z\subseteq \mathbb C^m$ of $z$.  
Note that we found  
holomorphic coordinates near each point of $Z_1$, so we have:

\begin{proclaim}{Corollary}
The union $Z_1$ of all codimension $\le 1$ strata, with the 
restriction of the sheaf $\mathfrak F_{V/G}$,  
is a complex manifold. \qed
\end{proclaim}
 
\subsection*{\nmb.{3.6}. The reflection divisor of a meromorphic 
tensor field on $Z_1$} 
Let $\Ga_{\mathcal M}(T^p_q(Z_1))$ be the space of meromorphic tensor 
fields (i.e.\  
with meromorphic coefficient functions in local holomorphic 
coordinates on the 
complex manifold $Z_1$), and let $P\in \Ga_{\mathcal M}(T^p_q(Z_1))$.

Let $S$ be an irreducible component of 
$Z_1\setminus Z_0$ and let $z\in S$. Local coordinates 
$y^1,\dots, y^n$ on $U\subset Z_1$, centered at $z$, are called 
adapted to the stratification of $Z_1$ if  $S=\{y^n=0\}$ near $z$. 
By definition the coordinates $y^1,\dots,y^n$ from \nmb!{3.5} 
have this property. Denote by $\mathcal O_z$ the ring of germs of 
holomorphic functions and by $\mathcal M_z$ the field of
germs of meromorphic functions, both at $z\in Z_1$.

Let $y^1,\dots,y^n$ be local coordinates on $U\subset Z_1$, centered 
at $z$, adapted to the stratification of $Z_1$ . 
Then on $U$ the meromorphic tensor field $P$ is given by 
$$
P|_U=\sum_{i_1,\dots,i_p,j_1,\dots,j_q} P^{i_1\dots i_p}_{j_1\dots j_q}
\frac{\partial}{\partial y^{i_1}}\otimes\dots\otimes
\frac{\partial}{\partial y^{i_p}}\otimes dy^{j_1}\otimes dy^{j_q}.
$$  
where the $P^{i_1\dots i_p}_{j_1\dots j_q}$ are meromorphic on $U$. 
Let us fix one nonzero summand of the right hand side: 
for the coefficient function we have 
$P^{i_1\dots i_p}_{j_1\dots j_q}=(y^n)^m f$ for some integer $m$ such 
that the germs at $z$ of $y^n$, $g$, and $h$ are pairwise relatively 
prime in $\mathcal O_z$ where $f=g/h\in \mathcal M_z$. 
Suppose that
the factor $\frac{\partial}{\partial y^n}$ appears exactly $p'$ times 
and the factor  
$dy^n$ appears exactly $q'$ times in this summand. The integer  
$$
\mu=mr+(q'-p')(r-1),
$$ 
a priori depending on $z$, is constant along an open dense subset of  
$S$ and it is called the 
\idx{\it reflection residuum} of the summand at $S$.  
Finally let $\mu_S(P)$ be the minimum of the reflection residua at $S$ 
of all summands of $P$ in the representation of $P$.

Let $\tilde y^1,\dots,\tilde y^n$ be arbitrary local coordinates on 
$U\subset Z_1$, centered at $z$, adapted to the stratification of $Z_1$. 
In a neighborhood of 
$z$ we have $y^n=f\tilde y^n$, where $f$ is a holomorphic function such that 
$f(z)\ne 0$. Remark that $\tilde y^n$ divides  
$\frac{\partial y^n}{\partial\tilde y^i}$ 
and $\frac{\partial\tilde y^n}{\partial y^i}$ $(i=1,\dots,n)$ in 
$\mathcal O_z$. 
A straightforward calculation using the above remark shows 
that the values of $\mu_S(P)$ calculated in the 
coordinates $\tilde y^i$ and in the coordinates $y^i$ are the same.  
Then $\mu_S(P)$ does not depend on the choice of the system of 
local coordinates adapted to the stratification of $Z_1$. For 
details see \cite{LMP}: there we checked this in the algebraic 
geometry setting where the use of tensor fields is less familiar.

We now can define the \idx{\it reflection divisor} 
$$
\on{div}_D(P)=\on{div}_{D_{V/G}}(P)\in \on{Div}(U)
$$
as follows: 
take the divisor $\on{div}(P)$, and for each irreducible 
component $S$ of $Z_1\setminus Z_0$ do the following:
if $S$ appears in the support of $\on{div}(P)\in \on{Div}(U)$, replace its 
coefficient by $\mu_S(P)$; if it does not appear, add 
$\mu_S(P).S$ to it. If $S$ is not contained in 
$Z_1\setminus Z_0$, we keep its coefficient in $\on{div}(P)$.

Finally we glue the global \idx{\it reflection divisor} 
$\on{div}_D(P)\in \on{Div}(Z_1)$ from the local ones, using a 
holomorphic atlas for $Z_1$. 

 
\begin{proclaim}{\nmb.{3.7}. Theorem}  
Let $G\subset GL(V)$ be a finite group, with 
reflection divisor $D=D_{V/G}=D_Z$. Then we have:
\begin{itemize}
\item 
       Let $P$ be a holomorphic $G$-invariant tensor field on 
       $V$. Then the reflection divisor $\on{div}_D(\pi_*P)\ge 0$.
\item 
       Let $Q\in \Ga_{\mathcal M}(T^p_q(Z_1))$ be a 
       meromorphic tensor 
       field on $Z_1$. Then the $G$-invariant meromorphic tensor field 
       $\pi^*Q$ extends to a holomorphic $G$-invariant tensor field 
       on $V$ if and only if 
       $\on{div}_D(Q)\ge 0$. 
\end{itemize}
The above remains true for $G$-invariant holomorphic tensor  
fields 
defined in a $G$-stable open subset of $V$. 
\end{proclaim} 
\begin{demo}{Proof} 
This follows directly from Hartogs' extension theorem,
the basic example \nmb!{3.4} using 
$y^1,\dots,y^{n-1}$ as dummy variables, and the definition of the 
reflection divisor $\on{div}_D(P)$ as explained in \nmb!{3.6}. 
\qed\end{demo} 
 
\begin{proclaim}{\nmb.{3.9}. Corollary} 
The mapping $\si$ establishes an injective correspondence between the space  
of holomorphic $G$-invariant tensor fields of type $\binom{p}q$ on 
$V$ which are 
skew-symmetric with respect to the covariant entries, and the space of 
holomorphic tensor fields on $Z_1$ of the same type  
and the same skew-symmetry condition. If $p=0$ the correspondence is 
bijective.  
 
The above remains true for $G$-invariant holomorphic tensor  
fields defined in a $G$-stable open subset of $V$. 
\end{proclaim} 
 
\begin{demo}{Proof} Let $P$ be a holomorphic $G$-invariant tensor field 
on $V$ satisfying the conditions of the corollary. For each nonzero 
decomposable summand of $\pi_*P$ take the integers 
$m$, $p'$, and $q'$ defined in \nmb!{3.6}. 
By skew symmetry of $P$ we have $q'\le 1$. 
By Theorem \nmb!{3.7} we get $\on{div}_D(\pi_*P)\ge 0$ and thus
$mr\ge (p'-q')(r-1)>-r$. So $m\ge 0$ and the summand is holomorphic 
on $Z_1$. 

If $Q$ is a holomorphic differential form on $Z_1$ its pullback 
$\si^*Q$ is a $G$-invariant holomorphic form on $\si^{-1}(Z_1)$ and 
then has a holomorphic extension to the whole of $V$. 
\qed\end{demo} 
 
\subsection*{\nmb.{3.10}. Remarks} 
Note that Corollary \nmb!{3.9} is a generalization of Solomon's   
theorem (see \cit!{6}): 
{\it 
If $G\subset GL(V)$ is a finite complex reflection group then 
every $G$-invariant polynomial exterior $q$-form 
$\om$ on $V$ can be written as $\om = \si^*\ph$ for a polynomial 
$q$-form $\ph$ on $\mathbb C^n$, where 
$\si=(\si^1,\dots,\si^n):V\to \mathbb C^n$ is the mapping consisting 
of a minimal system of homogeneous generators of $\mathbb C[V]^G$.} 

Actually, in  
the case of a reflection group $Z=\mathbb C^n$ and each holomorphic   
$\binom pq$-tensor field $Q$ on $Z_1$ has a holomorphic extension to  
$Z$ by Hartogs' extension theorem. 
 
\section*{\nmb0{4}. Invariant complex connections} 
 
\subsection*{\nmb.{4.1}} 
Let $\Gamma$ be a holomorphic $G$-invariant complex connection on $V$. 
Then the image $\sigma_*\Gamma$ of $\Gamma$ under the map $\sigma$ 
defines a holomorphic complex connection on $Z_0$. 
 
Let $z\in Z_1\setminus Z_0$, $v\in \sigma^{-1}(z)$, and $r$ the order of 
$G_v$. Consider the coordinates $z^i$ in $V$ defined in \nmb!{3.5}. 
Denote by $\Gamma^i_{jk}$ the components of the connection $\Gamma$ with 
respect to these coordinates. By assumption, the $\Gamma^i_{jk}$ are 
holomorphic functions on $V$. Recall the standard formula for the 
image $\gamma$ of $\Gamma$ under a holomorphic diffeomorphism 
$f=(y^a(x^i))$  
$$ 
\gamma^a_{bc}\circ f=\frac{\partial y^a}{\partial x^i}\frac{\partial 
x^j} 
{\partial y^b}\frac{\partial x^k}{\partial y^c}\Gamma^i_{jk}(x^l)- 
\frac{\partial^2y^a}{\partial x^i\partial x^j} 
\frac{\partial x^i}{\partial y^b}\frac{\partial x^j}{\partial 
y^c}. 
$$ 
Remark that the similar formula is true for the transformation 
of the components of connection under the change of coordinates. 

Consider the generator $g$ of the cyclic group $G_v$ given by 
\nmb!{3.5}. 
Since $g$ acts linearly, the connection reacts to it like a  
$\binom 12$-tensor field. Thus by the considerations of \nmb!{3.4} we  
get in the notation of \nmb!{3.5}, where $i,j,k=1,\dots,n-1$: 
\begin{gather*} 
\Gamma^i_{jk}=\tilde\Gamma^i_{jk}\o\sigma,\quad  
\Gamma^n_{jk}=\frac{1}{r}z^n\tilde\Gamma^n_{jk}\o\sigma,\quad 
\Gamma^i_{jn}=r(z^n)^{r-1}\tilde\Gamma^i_{jn}\o\sigma,\\ 
\Gamma^i_{nk}=r(z^n)^{r-1}\tilde\Gamma^i_{nk}\o\sigma,\quad 
\Gamma^n_{jn}=\tilde\Gamma^n_{jn}\o\sigma,\quad  
\Gamma^n_{nk}=\tilde\Gamma^n_{nk}\o\sigma,\\ 
\Gamma^i_{nn}=r^2(z^n)^{r-2}\tilde\Gamma^i_{nn}\o\sigma,\quad 
\Gamma^n_{nn}=r(z^n)^{r-1}\tilde\Gamma^n_{nn}\o\sigma, 
\end{gather*} 
where the $\tilde\Gamma^a_{bc}$ are holomorphic functions of the 
coordinates $y^a$ $(a=1,\dots,n)$ introduced in \nmb.{3.5}. 
 
Using the transformation formula for connections,  
we get the following formulas for the 
components $\gamma^a_{bc}$ of the meromorphic connection  
$\sigma_*\Gamma$ with respect 
to the coordinates $y^a$ 

{\begin{gather*}\tag{\nmb!{4.1}.1} 
\gamma^i_{jk}=\tilde\Gamma^i_{jk},\quad  
\gamma^n_{jk}=y^n\tilde\Gamma^n_{jk},\quad  
\gamma^i_{jn}=\tilde\Gamma^i_{jn},\quad  
\gamma^i_{nk}=\tilde\Gamma^i_{nk},\\ 
\gamma^n_{jn}=\tilde\Gamma^n_{jn},\quad  
\gamma^n_{nk}=\tilde\Gamma^n_{nk},\quad  
\gamma^i_{nn}=\frac1{y^n}\tilde\Gamma^i_{nn},\quad 
\gamma^n_{nn}=\tilde\Gamma^n_{nn}-\frac{r-1}{ry^n}. 
\end{gather*}}

Let $\tilde y^a$ for $a=1,\dots,n$ be other local coordinates 
centered at $z$ and adapted to the stratification of $Z_1$. 
Then in a neighborhood of $z$ we have 
$$
y^n=f\tilde y^n,\quad \tilde y^n=\tilde fy^n,
$$
where $f$ and $\tilde f$ are holomorphic functions in a neighborhood of 
$z$ and $\tilde ff=1$. Then we have 
$$
\frac{\partial y^n}{\partial\tilde y^i}=\frac{\partial f}{\partial\tilde y^i}
\tilde y^n,
\quad \frac{\partial\tilde y^n}{\partial y^i}=\frac{\partial f}{\partial y^i}
y^n\quad (i=1,\dots,n-1)  
$$
and on $S=\{y^n=0\}$
$$
\frac{\partial y^n}{\partial\tilde y^n}=f,
\quad\frac{\tilde \partial y^n}{\partial y^n}=\tilde f.
$$   
Using these formulas one can check 
that in the coordinates $\tilde y^a$ the formulas \nmb!{4.1}.1  have the 
same form as in the coordinates $y^a$. 
For example, for the new component $\tilde\gamma^n_{nn}$ we have 
$$
\tilde\gamma^n_{nn}+\frac{r-1}{r\tilde y^n}=
\frac{(r-1)\biggl(1-\tilde f\frac{\partial\tilde y^n}{\partial y^n}
\biggl(\frac{\partial y^n}
{\partial\tilde y^n}\biggr)^2\biggr)}{ry^n\tilde f}+h,
$$
where $h$ is a holomorphic function near $z$.
Since on $S=\{y^n=0\}$ we have 
$$
1-\tilde f\frac{\partial\tilde y^n}{\partial y^n}\biggl(\frac{\partial y^n}
{\partial\tilde y^n}\biggr)^2=1-\tilde f^2f^2=0,
$$ 
$y^n$ divides in $\mathcal O_z$ the function  
$$
1-\tilde f\frac{\partial\tilde y^n}{\partial y^n}\biggl(\frac{\partial y^n}
{\partial\tilde y^n}\biggr)^2.
$$
Thus  
$$
\tilde\gamma^n_{nn}+\frac{r-1}{r\tilde y^n}
$$ 
is holomorphic in a neighborhood of $z$.

\begin{proclaim}{\nmb.{4.2}. Theorem} 
Let $\gamma$ be a holomorphic complex 
linear connection on $Z_0$ such that for each $z\in Z_1\setminus Z_0$ it 
has an extension to a neighborhood of $z$ whose components in the 
coordinates adapted to the stratification of $Z_1$ 
are defined by the formulas \nmb!{4.1}.1  
where $\tilde\Gamma^a_{bc}$ are holomorphic. 
Then there is a unique $G$-invariant holomorphic complex linear 
connection $\Gamma$ on $V$ such that $\sigma_*\Gamma$ coincides 
with $\gamma$ on $Z_0$. This remains true if we replace $V$ by a 
$G$-open subset of $G$. 
\end{proclaim} 
 
\begin{demo}{Proof} 
Since $\si$ is étale on the principal stratum, there is a unique 
$G$-invariant complex linear connection $\Gamma_0$ on $\sigma^{-1}(Z_0)$ 
such that $\sigma_*\Gamma_0=\gamma$.  The condition of the theorem implies 
that the connection $\Gamma_0$ has a holomorphic extension to 
$\sigma^{-1}(Z_1)$. 
Then by Hartogs' extension theorem the connection $\Gamma_0$ has a 
unique holomorphic extension $\Gamma$ to the whole of $V$. 
\qed\end{demo} 
 
\section*{\nmb0{5}. Lifts of diffeomorphisms of orbit spaces} 
 
\subsection*{\nmb.{5.1}} 
Let $G$ and $G'$ be finite subgroups of $GL(V)$ and $GL(V')$ and let 
$F$ be a holomorphic diffeomorphism $V\to V'$ which maps $G$-orbits to 
$G'$-orbits bijectively.  
Then the map $F$ induces an isomorphism $f$ of the sheaves 
$\mathfrak F_{V/G}\to\mathfrak F_{V'/G'}$, 
i.e.\ a holomorphic diffeomorphism of orbit spaces $V/G$ and $V'/G'$. 

\begin{proclaim}{Lemma}  
There is a unique isomorphism $a:G\to G'$ such that  
$F\o g=a(g)\o F$ for every $g\in G$.  
\end{proclaim}

Note that $a$ and its inverse $a^{-1}$ map complex reflections to 
complex reflections.  

\begin{demo}{Proof} 
The cardinalities of the two groups are the same since $F$ maps a generic 
regular orbit to a regular orbit. Consequently, it maps regular points to 
regular points and we have $\sigma'\o F= 
f\o \sigma:V\to V'/G'$ for a holomorphic diffeomorphism  
$f:V/G\to V'/G'$, where 
$\sigma:V\to V/G$ and $\sigma':V'\to V'/G'$ are the quotient projections.  

Fix some $G$-regular $v\in V$. Then $F(v)$ and 
$F(gv)$ for $g\in G$ are regular points of $V'$ of the same orbit. 
Therefore, there is a unique   
$a(g)\in G$ such that $F(gv)=a(g)(F(v))$.
We have 
$\sigma'\o F\o g = f\o \sigma\o g =f\o \sigma = \sigma'\o F = \sigma'\o a(g)\o F$.
Since $\sigma'$ is étale on $V'_{\text{reg}}$ we see that 
$F\o g= a(g)\o F$ locally near $v$ and thus globally. 
By uniqueness, the map 
$g\to a(g)$ is an isomorphism of $G$ onto $G'$. 
\qed\end{demo} 
 
In this section we study when a diffeomorphism $f$ of the orbit spaces 
$Z\to Z'$ has a holomorphic lift $F$. 
 
\begin{proclaim}{\nmb.{5.2}. Corollary}  
Let $F:V\to V$ be a holomorphic diffeomorphism  
which maps $G$-orbits onto $G'$-orbits,  
and $f:Z\to Z'$ the corresponding holomorphic diffeomorphism of  
the orbit spaces. Then $f$ maps the isotropy type stratification 
of $Z$ onto that of $Z'$ and, moreover, it maps $D_Z$ to $D_{Z'}$.  
\end{proclaim} 

\begin{demo}{Proof}  
This follows from Lemma \nmb!{5.1} and the definition 
\nmb!{3.3} of the reflection divisor.
\qed\end{demo} 
 
\begin{proclaim}{\nmb.{5.3}. Theorem}  
Let $G$ and $G'$ be two finite subgroups of $GL(V)$ and let 
$f:Z\to Z'$ be a holomorphic diffeomorphism of the corresponding 
orbit spaces such that $f(Z_0)=Z_0'$ and $f_*(D_Z)=D_{Z'}$. 
If $Q$ is a holomorphic tensor field of 
type $\binom{p}q$ on $Z_0$ which satisfies the conditions of Theorem  
\nmb!{3.7},  
then $f_*(Q)$ also satisfies these conditions on $Z_0'$ and thus 
there exists a unique $G'$-invariant holomorphic tensor field $Q'$ of  
type $\binom{p}{q}$ such that $\sigma'_*Q'$ coincides with $f_*Q$ on $Z_0'$. 
 
This is also true for holomorphic connections if we replace Theorem \nmb!{3.7} 
by Theorem \nmb!{4.2}. 
The theorem remains true if we replace $V$ by 
invariant open subsets of $V$.
\end{proclaim} 
 
\begin{demo}{Proof}  
Since $f(Z_0)=Z_0'$ the tensor field $f_*Q$ is also holomorphic 
on $Z_0'$. 
Let $z\in Z_1\setminus Z_0$. Then there is a complex space   
$S\in R_Z$ of codimension $1$ such that $z\in S$.  
By assumption $f(z)\in Z_1'\setminus Z_0'$ and $f(z)\in f(S)\in R_{Z'}$ and 
$r_z=e_S=e_{f(S)}=r_{f(z)}$. Now, obviously $f_*Q$ satisfies the conditions of 
Theorem \nmb!{3.7} at $f(x)$. Thus there exists a $G'$-invariant holomorphic  
tensor field $Q'$ on $V$ with $\si'_*Q'=f_*Q$. 
 
A similar argument applies to connections.
\qed\end{demo} 
 
\begin{proclaim}{\nmb.{5.4} Theorem} 
Let $G$ and $G'$ be two finite subgroups of $GL(V)$.  
Let $f:Z\to Z'$ be a holomorphic diffeomorphism of  
the orbit spaces such that $f(Z_0)=Z_0'$ and  
$f_*(D_Z)=D_{Z'}$. 
 
Then $f$ lifts to a holomorphic diffeomorphism $F:V\to V$, 
i.e.\ $\sigma'\o F=f\o\sigma$. 
 
The local version is also true. Namely, 
if $B$ is a ball in the vector space $V$ centered at $0$ (for an 
invariant Hermitian metric),  
$U=\sigma(B)$, 
and $f:U\to Z'$ is a local holomorphic diffeomorphism of $U$ onto 
a neighborhood $U'$ of $\sigma'(0)$ such that 
$f(U\cap Z_0)=U'\cap Z_0'$ and $f$ maps $D_Z\cap U$ to $D_{Z'}\cap U'$, 
then there is a holomorphic lift $F:B\to V$. 
\end{proclaim} 
 
\begin{demo}{Proof} 
Let $\Gamma$ be the natural flat connection on $V$. Then $\Gamma$ is 
uniquely defined by the holomorphic connection $\sigma_*\Gamma$ on $Z_0$ 
which satisfies the conditions of Theorem \nmb!{4.2}. 
By Theorem \nmb!{5.3} there is a unique $G$-invariant holomorphic complex 
linear connection $\Ga'$ on $V$ such that $\sigma'_*\Ga'$ 
coincides with $f_*(\sigma_*\Gamma)$ on $Z_0'$. 
It is evident that $\Ga'$ is a torsion free flat connection, since 
$\Ga$ is it and $\Ga'$ is locally isomorphic to $\Ga$ on an open 
dense subset.  
 
Let $v\in V$ be $G$-regular and let $v'\in V$ be $G'$-regular, such that 
$(f\o\sigma)(v)=\sigma'(v')$. Then there is a biholomorphic 
map $F$ of a neighborhood $W$ of $v$ onto a neighborhood of $v'$ such  
that 
$\sigma'\o F=f\o\sigma$ on $W$ and 
$F(v)=v'$. Moreover by construction $F$ is a locally affine map of the 
affine space $(V,\Ga)$ into $(V,\Ga')$ equipped with the above structures of 
locally affine spaces, thus we have  
\begin{equation}\label{exponential} 
F = \exp^{\Ga'}_{v'}\o T_vF \o (\exp^{\Ga}_v)\i 
\end{equation} 
where $\on{exp}^{\Ga}_v:T_vV \to V$ is the holomorphic geodesic exponential 
mapping centered at $v$ given by the connection $\Ga$ and its induced spray.  
It is globally defined, thus complete and a holomorphic diffeomorphism  
since $\Ga$ is the standard flat connection. Likewise $\exp^{\Ga'}_{v'}$ is 
the holomorphic exponential mapping of the flat connection $\Ga'$. The 
formula above extends $F$ to a globally defined holomorphic  
mapping  
if $\exp^{\Ga'}_{v'}:T_vV\to V$ is also globally defined (complete). Assume 
for contradiction that this is not the case. Let $F$ be maximally extended 
by equation (\ref{exponential}); it still projects to $f:Z\to Z'$.  
We consider  
$\exp^{\Ga'}_{v'}$ as a real exponential mapping, and then there is a 
real geodesic which reaches infinity in finite time and this is the image 
under $F$ of a finite part $\exp^{\Ga}_v([0,t_0)w)$ of a real geodesic of  
$\Ga$ emanating at $v$. The sequence $\exp^{\Ga}_v((t_0-1/n)w)$  
converges to $\exp^{\Ga}_v(t_0w)$ in $V$, but its image under $F$  
diverges to infinity by assumption.  
On the other hand, the image under $F$ is contained in the set  
$(\si')\i(f\si(\exp^{\Ga}_v([0,t_0]w)))$ which is compact since  
$\si'$ is a proper mapping. Contradiction.  
 
Any holomorphic lift $F$ of a holomorphic diffeomorphism $f$ is a  
holomorphic   
diffeomorphism 
of $V$ which maps $G$-orbits onto $G'$ orbits, by the following 
argument: 
Let $F'$ be a holomorphic lift of $f^{-1}$. 
Evidently the map $F'\o F$ preserves each $G$-orbit. 
Then, for a $G$-regular point $v\in V$, there is a 
$g\in G$ such 
that $F'\o F=g$ in a neighborhood of $v$ and,  
then, 
on the whole of $V$.  
Similarly $F\o F'=g'\in G'$. 
This implies that $F$ is a holomorphic 
diffeomorphism of $V$. By definition the lift $F$ respects the  
partitions of 
$V$ into orbits. 
\qed\end{demo} 

We give a second proof of Theorem \nmb!{5.4} based on the 
known results about the fundamental groups of $V_{\text{reg}}$ and 
$Z_0$ for finite complex reflection groups. It is an extension of the 
proof of \cite{Ly}, using results of \cite{1}.
 
\begin{proclaim}{\nmb.{5.5}. Lemma}  
Let $G$ and $G'$ be two finite subgroups of $GL(V)$ and let $f:Z\to Z'$ 
be a holomorphic diffeomorphism of the corresponding orbit spaces. 
Suppose $v_0\in V_{\text{reg}}$, $v_0'\in V'_{\text{reg}}$, and  
$f\o\sigma(v_0)=\sigma'(v_0')$.   
If the image of the fundamental group 
$\pi_1(V_{\text{reg}},v_0)$ under $f\o\sigma$ is contained in the subgroup 
$\sigma'_*(\pi_1(V_{\text{reg}}),v_0')$ of $\pi_1(Z'_0,\sigma'(v'_0))$, the 
holomorphic lift of $f\o\sigma$ mapping $v_0$ to $v_0'$ exists. 
\end{proclaim} 
\begin{demo}{Proof}  
Consider the restriction $\varphi$ of the map  
$f\o\sigma$ to $V_{\text{reg}}$. 
Since the restriction of $\sigma$ to $V_{\text{reg}}$ is a covering 
map onto $Z_0$, 
the condition of the lemma implies that there is a holomorphic 
lift $F_0$ of the map $\varphi$ to $V_{\text{reg}}$. 
The map $F_0$ is bounded on $B\cap V_{\text{reg}}$ for each compact 
ball $B$ in $V$ since its image is contained in the compact set 
$(\si')^{-1}(f(\si(B)))$. 
Then by the Riemann extension theorem $F_0$ has a holomorphic 
extension $F$ to $V$ which is the required holomorphic lift of $f$. 
\qed\end{demo} 
 
\subsection*{\nmb.{5.6}} 
Next we prove Theorem \nmb!{5.4} in the case when the group $G$ is generated  
by 
complex reflections. Put 
\begin{equation*} 
B:=\pi_1(Z_0)\quad\text{and}\quad P:=\pi_1(V_{\text{reg}}). 
\end{equation*} 
The groups $B$ and $P$ are called the {\it braid group} and the {\it pure 
braid group} associated to $G$, respectively. It is clear that the map 
$\sigma$ induces an isomorphism of $P$ onto a subgroup of $B$. 
 
The following results about the groups $B$ and $P$ are well known 
(see, for example, \cit!{1}). The braid group $B$ is generated by 
those elements which are represented by loops around the hypersurfaces 
$\sigma(H)$ for $H\in\mathfrak H$. The pure braid group $P$ is generated 
by the elements of $B$ of the type $s^{e_H}$, where $s$ is any of 
the above generators of $B$ represented by a loop around the  
hypersurface $\sigma(H)$. 
This implies the following 
 
\begin{proclaim}{Proposition}  
Suppose the group $G$ is generated by complex reflections.  
Let $f$ be a holomorphic 
diffeomorphism of the orbit space $Z=\mathbb C^n$ with  
$f(Z_0)=Z_0$ which also preserves $D_Z$. 
Then $f|_{Z_0}$ preserves the subgroup $P$ of $B$.\qed 
\end{proclaim} 
 
The following proposition is an immediate consequence of Lemma  
\nmb!{5.5} and Proposition \nmb!{5.6}. 
 
\begin{proclaim}{\nmb.{5.7}. Proposition} 
Suppose the groups $G$ and $G'$ are  
generated by complex reflections. Let $f:Z\to Z'$ be a holomorphic  
diffeomorphism between the corresponding orbit spaces, such that $f(Z_0)=Z'_0$ 
and $f_*(D_Z)=D_{Z'}$. 
 
Then $f$ has a holomorphic lift $F$ to $V$.\qed 
\end{proclaim} 
 
\begin{demo}{Second proof of \nmb!{5.4}} 
Now let $G\subset GL(V)$ be a finite group and let $G_1$ be the 
subgroup generated by all  
complex reflections in $G$. Clearly $G_1$ is a normal subgroup  
of $G$. Let $G_2=G/G_1$. 
Let $\sigma_1^1,\dots,\sigma_1^n$ be a system of homogeneous generators 
of $\mathbb C[V]^{G_1}$ and $\sigma_1:V\to\mathbb C^n$ the corresponding 
orbit map. Then the action of $G$ on $V$ induces the action of the group 
$G_2$ on $V_1:=\mathbb C^n=\sigma_1(V)$. Since each representation of the  
group 
$G_2$ is completely reducible, by standard arguments of invariant  
theory, we may assume that the generators 
$\sigma_1^i$'s are chosen in such a way that the above action of 
$G_2$ on $V_1=\mathbb C^n$  
is linear. Then the representation of $G_2$ on $V_1$  
contains no complex reflections. 
Let $\sigma_2^1,\dots,\sigma_2^m$ be a system of homogeneous generators 
of $\mathbb C[V_1]^{G_2}$ and $\sigma_2:V_1\to\mathbb C^m$ the corresponding 
orbit map. Then $\sigma^i=\sigma_2^i\o\sigma_1$ 
$(i=1,\dots,m)$ is a system of generators of $\mathbb C[V]^G$ with 
orbit map $\sigma=\sigma_2\o\sigma_1$. 
Similarly for $G'$. 
 
Let $f:Z\to Z'$ be a holomorphic diffeomorphism, such that  
$f(Z_0)=Z'_0$ and $f_*(D_Z)=D_{Z'}$. 
Since the group 
$G_2$ contains no complex reflections the set $V_{1,\text{reg}}$ of regular 
points of the action of $G_2$ on $V_1$ is obtained from $V_1$ by  
removing 
some subsets of codimension $\geq 2$.  
And similarly for $G'$.  
Then the fundamental group 
$\pi_1(V_{1,\text{reg}})=\pi_1(V_1)=0$ is trivial and by lemma 
\nmb!{5.5} the diffeomorphism 
$f$ has a holomorphic lift $F_1:V_1\to V'_1$ which is a holomorphic  
diffeomorphism mapping the principal stratum to the principal 
stratum, and the reflection divisor to the reflection divisor,
since $G_2$ contains no complex reflections on $V_1$.
Thus the diffeomorphism $F_1$ has a holomorphic lift to $V$ by 
Proposition \nmb!{5.7}, which is a holomorphic lift of $f$. 
\qed\end{demo}

\section*{ 
\nmb0{6}. An intrinsic characterization of a complex orbifold} 
 
We recall the definition of orbifold. 
\begin{proclaim}{\nmb.{6.1}. Definition} 
\cit!{7} 
Let $X$ be a Hausdorff space. An atlas of a smooth $n$-dimensional  
orbifold 
on $X$ is a family $\{U_i\}_{i\in I}$ of open sets that satisfy: 
\begin{enumerate} 
\item $\{U_i\}_{i\in I}$ is an open cover of X. 
\item For each $i\in I$ we have a local uniformizing system consisting of 
a triple $(\tilde U_i,G_i,\varphi_i)$, where $\tilde U_i$ is a 
connected open subset of $\mathbb R^n$ containing the origin, 
$G_i$ is a finite 
group of diffeomorphisms acting effectively and properly on $\tilde  
U_i$, and 
$\varphi_i:\tilde U_i\to U_i$ is a continuous map of $\tilde U_i$ onto  
$U_i$ 
such that $\varphi_i\o g=\varphi_i$ for all $g\in G_i$ and 
the induced map of $\tilde U_i/G_i$ onto $U_i$ is a homeomorphism. 
The finite group $G_i$ is called a local uniformizing group. 
\item Given $\tilde x_i\in\tilde U_i$ and $\tilde x_j\in\tilde U_j$ such  
that 
$\varphi_i(\tilde x_i)=\varphi_j(\tilde x_j)$, there is a  
diffeomorphism 
$g_{ij}:\tilde V_j\to\tilde V_i$ from a neighborhood 
$\tilde V_j\subseteq\tilde U_j$ of $\tilde x_j$ onto a neighborhood 
$\tilde V_i\subseteq\tilde U_i$ of $\tilde x_i$ such that 
$\varphi_j=\varphi_i\o g_{ij}$. 
\end{enumerate} 
Two atlases are equivalent if their union is again an atlas of a smooth 
orbifold on $X$. An orbifold is the space $X$ with an equivalence class 
of atlases of smooth orbifolds on $X$. 
\end{proclaim} 
If we take in the definition of orbifold $\mathbb C^n$ instead of 
$\mathbb R^n$ and require that $G_i$ is a finite group of holomorphic 
diffeomorphisms acting effectively and properly on $\tilde U_i$ and the 
maps $g_{ij}$ are biholomorphic, we get the definition of complex 
analytic $n$-dimensional orbifold. 
 
\begin{proclaim}{\nmb.{6.2}. Theorem} \cit!{7} Let $M$ be a smooth manifold and $G$ a 
proper discontinuous group of diffeomorphisms of $M$. Then the orbit  
space 
$M/G$ has a natural structure of smooth $n$-dimensional orbifold. 
If $M$ is a complex $n$-dimensional manifold 
and $G$ is a group of holomorphic diffeomorphisms of $M$, the orbit  
space 
$M/G$ is a complex $n$-dimensional orbifold. 
\end{proclaim} 
 
\subsection*{\nmb.{6.3} Definitions} 
In the definition of atlas of a complex orbifold on $X$ we can 
always take $\tilde U_i$ to be balls of the space $\mathbb C^n$ (with 
respect to some Hermitian metric) centered 
at the origin and the finite subgroups $G_i$ to be subgroups of the 
$GL(n)$ acting naturally on $\mathbb C^n$. In the sequel we 
consider atlases of complex orbifolds satisfying these  
conditions. 
 
Let $X$ be a complex orbifold with an atlas 
$(\tilde U_i,G_i,\varphi_i)$. A function $f:U_i\to\mathbb C$ is called 
holomorphic if $f\o\varphi_i$ is a holomorphic function on $\tilde  
U_i$. 
The germs of holomorphic functions on $X$ define a {\it sheaf} $\mathfrak F_X$ on  
$X$. 
It is evident that the sheaf $\mathfrak F_X$ depends only on the structure 
of complex orbifold on $X$. 
 
Consider a uniformizing system $(\tilde U_i,G_i,\varphi_i)$ of the 
above atlas and the corresponding action of $G_i$ on $\mathbb C^n$. Then 
we have the isotropy type stratification of the orbit space 
$\mathbb C^n/G_i$, the induced stratification of $U_i$, and the  
divisor $D_{U_i}$. 
 
By corollary \nmb!{5.2} we get the {\it stratification on $X$} by  
gluing the strata  
on the $U_i$'s. Denote by $X_0$ the principal stratum of this  
stratification. 
By definition, for each $x\in X_0$, for each uniformizing system 
$(\tilde U_i,G_i,\varphi_i)$, and for each  
$y\in\tilde U_i$ such that 
$\varphi_i(y)=x$, the isotropy group $G_y$ of $y$ is trivial. Note  
that $X_0$ is a complex manifold. 
Note that $X_1$ is also a complex manifold 
since this holds locally as noted in \nmb!{3.5}. 
 
Denote by $R_X$ the {\it set of all strata of codimension $1$} of $X$.  
Since the pullbacks of the reflection divisors $D_{U_i}$ to $U_i\cap U_j$ agree 
by \nmb!{5.2} we may glue them into the reflection divisor $D_X$ on $X_1$.  
 
\begin{proclaim}{\nmb.{6.4}. Definition} 
Let $X$ and $\tilde X$ be two smooth  
orbifolds. 
The orbifold $\tilde X$ is called a covering orbifold for $X$ 
with a projection $p:\tilde X\to X$ if $p$ is a continuous map of  
underlying 
topological spaces and  each point $x\in X$ has a neighborhood 
$U=\tilde U/G$ (where $\tilde U$ is an open subset of $\mathbb R^n$) 
for which each component $V_i$ of $p^{-1}(U)$ is isomorphic to 
$\tilde U/G_i$, where $G_i\subseteq G$ is some subgroup. The 
above isomorphisms $U=\tilde U/G$ and $V_i=\tilde U/G_i$ must 
respect the projections. 
\end{proclaim} 
Note that the projection $p$ in the above definition is not 
necessarily a covering of the 
underlying topological spaces. 
It is clear that a covering orbifold for a complex orbifold is  
a complex orbifold. 
Hereafter we suppose that all orbifolds and their covering orbifolds 
are connected. 
 
\begin{proclaim}{\nmb.{6.5}. Theorem} \cit!{7}  
An orbifold $X$ has a universal covering 
orbifold $p:\tilde X\to X$. More precisely, if $x\in X_0$, 
$\tilde x\in \tilde X_0$ and $p(\tilde x)=x$, 
for any other covering orbifold $p':\tilde X'\to X$ and 
$\tilde x'\in\tilde X'$ such that $p'(\tilde x')=x$ there is a cover 
$q:\tilde X\to\tilde X'$ such that $p=p'\o q$ and $q(\tilde  
x)=\tilde x'$. 
For any points $\tilde x,\tilde x'\in p^{-1}(x)$ there is a deck 
transformation of $\tilde X$ taking $\tilde x$ to $\tilde x'$. 
\end{proclaim} 
 
Now we prove the main theorem of this section. 
 
\begin{proclaim}{\nmb.{6.6}. Theorem} 
An $n$-dimensional complex orbifold $X$ is uniquely determined by  
the sheaf of holomorphic functions $\mathfrak F_X$, 
the principal stratum $X_0$, and the reflection divisor $D_X$. 
\end{proclaim} 
 
\begin{demo}{Proof}  For each $x\in X$, there exists $V=\mathbb C^m$, 
a finite group $G\subset GL(m)$,
a ball $B$ in $V$ centered at 0, an open subset $U$ of $X$ containing 
$x$, and an isomorphism $\psi:\pi(B)\to U$ between   the sheaves 
$\mathfrak F_Z|_{\pi(B)}$ and $\mathfrak F_X|_U$.
Consider the  map $\pi:V\to Z=V/G$, the stratum $Z_0$ and the 
reflection divisor $D_Z$. We suppose also that 
$\psi(Z_0\cap B/G)\subseteq X_0$ and $\ps_*(D_{\pi(B)})=D_U$. It 
suffices to prove that the germ of the uniformizing system 
$\{B,G,\psi\o\pi|B\}$ at $x$ is the germ of some uniformizing system 
of the orbifold $X$. 
 
Let $y\in V_{\text{reg}}\cap B$. Then the ring $\mathfrak F_Z(\pi(y))$ of 
germs of $\mathfrak F_Z$ at $\pi(y)$ is isomorphic to the ring of germs of 
holomorphic functions on $\mathbb C^n$ at $0$ and thus we have $m=n$. 
 
Consider the uniformizing system $(\tilde U_i,G_i,\ph_i)$ of the 
orbifold $X$, where $\tilde U_i$ is a ball in $\mathbb C^n$ centered 
at the origin, $G_i$ is a finite subgroup of the group $GL(n)$ acting 
naturally on $V=\mathbb C^n$, and where $\ph_i(0)=x$. Consider the 
map $\pi_i:V\to V/G_i$ given by some system of generators of 
$\mathbb C[V]^{G_i}$. We may assume that 
$\ph_i=\psi_i\o\pi_i|_{\tilde U_i}$, where 
$\ps_i:\mathfrak F_{\tilde U_i/G_i}\to\mathfrak F_{U_i}$ is an 
isomorphism of sheaves.
$$ \xymatrix{ 
 \mathbb C^n &B\ar@{->}[1,0]_{\pi} \ar@{_(->}[0,-1] \ar@{.>}[0,2]_{F}  &  & 
\tilde U_i\ar@{->}[1,0]^{\pi_i} \ar@{^(->}[0,1]  &\mathbb C^n \\ 
  &B/G\ar@{->}[1,1]_{\ps} \ar@{.>}[0,2]_{f}  &  &
  \tilde U_i/G_i\ar@{->}[1,-1]^{\ps_i}  &  \\ 
  &  &U &  &  \\ 
} 
$$ 
Then the maps $\ps$ and $\ps_i$ define a map (germ) $f$ of a holomorphic 
diffeomorphism $B/G$ to $U_i/G_i$ at $0:=\pi(0)$ such that 
$f(0)=0:=\pi_i(0)$. Then $f$ induces an 
isomorphism $\mathfrak F_{V/G}(0)\to\mathfrak F_{V/G_i}(0)$, it maps $(B/G)_0$ to 
$(\tilde U_i/G_i)_0$ and $f_*(D_{B/G})=D_{\tilde U_i/G_i}$. 
Thus by theorem \nmb!{5.4} there is a germ of a holomorphic diffeomorphism 
$F:B\to \tilde U_i$ which is equivariant for a suitable isomorphism $G\to G_i$. 
\qed\end{demo} 
 
\begin{proclaim}{\nmb.{6.7}. Corollary} Let $M$ be a complex simply connected 
manifold, $G$ a proper discontinuous group of holomorphic 
diffeomorphisms of $M$, and $\mathfrak F_X$ the corresponding sheaf 
on the orbifold $X=M/G$. The $G$-manifold $M$ is a universal covering 
orbifold for the orbifold $X$ and it is defined uniquely up to 
a natural isomorphism of universal coverings by the sheaf $\mathfrak F_X$, 
the principal stratum $X_0$, and by the reflection divisor $D_X$. 
\end{proclaim} 
\begin{demo}{Proof} Evidently the manifold $M$ is a covering orbifold for $X$. 
If $\tilde X$ is a universal covering orbifold for $X$, by definition  
\nmb!{6.4} 
there is a cover $q:\tilde X\to M$. By definition $\tilde X$ should be 
a manifold and $q$ a cover of manifolds. Therefore, $q$ is a  
diffeomorphism. 
Then the statement of the corollary follows from theorem \nmb!{6.6}. 
\qed\end{demo} 
An automorphism of the sheaf $\mathfrak F_X$ is called a holomorphic 
diffeomorphism of the orbit space $X$. 
Theorem \nmb!{6.5} and corollary \nmb!{6.7} imply the following analogue of 
Theorem \nmb!{5.4}. 
 
\begin{proclaim}{\nmb.{6.8}. Theorem} Let $M$ be a complex  
simply connected 
manifold, $G$ a proper discontinuous group of holomorphic 
diffeomorphisms of $M$, and $\mathfrak F_X$ the corresponding sheaf 
on the orbifold $X=M/G$. Each holomorphic diffeomorphism $f$ of the  
orbit space $X$ preserving $X_0$ and $D_X$ 
has a holomorphic lift $F$ to $M$, which is $G$-equivariant with respect 
to an automorphism of $G$. The lift $F$ is unique up to composition 
by an element of $G$.  
\end{proclaim} 
 
\begin{demo}{Proof} By theorem \nmb!{6.6} and corollary \nmb!{6.7} the  
manifold $M$ with the map 
$f\o p:M\to X$, where $p:M\to X$ is the projection, is a universal 
covering orbifold for $X$. Then there is a holomorphic diffeomorphism 
$F:M\to M$ such that $p\o F=f\o p$. The equivariance property 
holds locally by \nmb!{5.1}, thus globally. The lift is uniquely 
given by choosing $F(x)$ for a regular point $x$ in the orbit 
$f(p(x))$. 
\qed\end{demo} 
 
\subsection*{\nmb.{6.9}} 
Let $V$ be a complex vector space with a linear action of a finite  
group $G$. The group $\mathbb C^*$ acts on $V$ by homotheties and  
induces an action on $Z=V/G$. 
 
\begin{proclaim}{Corollary} 
In this situation, the $G$-module $V$ is uniquely defined up to a  
linear isomorphism by the sheaf $\mathfrak F_{V/G}$ with the action of  
$\mathbb C^*$, by $Z_0$, and the reflection divisor $D_Z$. \qed 
\end{proclaim}  

\begin{demo}{Proof}
Consider the orbit space $Z=V/G$ of a $G$-module $V$ with the sheaf 
$\mathfrak F_{V/G}$, regular stratum $Z_0$, reflection divisor $D_Z$, 
and the action of $\mathbb C^*$ induced by the action of 
$\mathbb C^*$ on $V$ by homotheties. Suppose  
that we have another $G'$-module $V'$ with the same data on 
$Z'=V'/G'$ such that there is a biholomorphic map $f:Z\to Z'$ 
preserving these data. By Theorem \nmb!{4.5} 
there is a biholomorphic lift 
$F: V\to V'$, and by lemma \nmb!{5.1} there is an isomorphism $a:G\to G'$ 
such that $F\circ g=a(g)\circ F$. Thus we may assume that  
$G=G'$, $V=V'$, $Z=Z'$, and $a$ is the identity map. 
By definition the pullback $A$ of the vector field on the orbit space 
$V/G$ defined by the action of the group $\mathbb C^*$ on $V/G$ 
coincides with the vector field on $V$ defined by the above action of 
the group $\mathbb C^*$ on $V$. By construction $F^*A=A$ and then the 
map $F$ commutes with the action of $\mathbb C^*$ on $V$, i.e. for 
each $t\in\mathbb C^*$ and $v\in V$ we have $F(tv)=tF(v)$. Since $F$ is 
biholomorphic it is a linear automorphism of the vector space 
$V$. By definition $F$ is then an automorphism of the $G$-module $V$. 
\qed\end{demo}
 
\subsection*{\nmb.{6.10}. Tensor fields and connections on orbifolds} 
The local results in section \nmb!{3} show that the correct 
definition of a 
$\binom{p}{q}$-tensor field $Q$ on an orbifold $X$ is as follows:  
$Q$ is a meromorphic 
$\binom{p}{q}$-tensor field on $X_1$ such that 
$\operatorname{div}_{D_X}(Q)\ge 0$. 
 
Likewise, we can define connections on orbifolds by requiring the local 
conditions of section \nmb!{4}.


\begin{thebibliography}{9} 


 \bibitem{Bier}
 E. Bierstone,
 \emph{ Lifting isotopies from orbit spaces}, 
 Topology 14 (1975), 245--252.
 
\bibitem{1} 
M. Brou\'e, G. Malle, R. Rouquier, 
\emph{ Complex reflection groups, braid groups, Hecke algebras}, 
J. reine angew. Math 
500 (1998) 
127-190.
 
\bibitem{2a} 
M. Losik, 
\emph{ Lifts of diffeomorphisms of orbit spaces for 
representations of compact Lie groups}, 
Geom. Dedicata, 88 (2001), 21-36. 

\bibitem{LMP} 
M. Losik, P.W. Michor, V.L. Popov,
\emph{ Invariant tensor fields and orbit varieties for 
finite algebraic transformation groups}, 
arXiv:math.AG/0206008. 

 \bibitem{2} 
 D. Luna, 
 \emph{ Slices \'etales}, 
 Bull. Soc. Math. France,  M\'emoire 33 (1973), 81-105. 
 
 \bibitem{3} 
 D. Luna, 
 \emph{ Sur certaines op\`erations diff\'erentiables des groupes de Lie}, 
 Amer. J. Math. 97 (1975), 172-181. 
 
  
 \bibitem{4} 
 D. Luna, 
 \emph{ Fonctions diff\'erentiables invariantes sous l'op\'eration d'une 
groupe r\'eductif}, 
 Ann. Inst. Fourier 
 26 (1976), 33-49. 

 \bibitem{Ly}
O. V. Lyashko,
\emph{ Geometry of bifurcation diagrams},
J. Soviet Math. 27 (1984), 2736-2759.
 
\bibitem{5} 
G.W. Schwarz, 
\emph{ Lifting smooth homotopies of orbit spaces}, 
Publ. Math. IHES 51 (1980), 37-136. 
 
  
\bibitem{6} 
L. Solomon, 
\emph{ Invariants of finite reflection groups}, 
Nagoya J. Math. 22 (1963), 57-64. 
 
  
 \bibitem{7} 
 W.P. Thurston, 
 \emph{ The geometry and topology of tree-manifolds}, 
 Lect. Notes 
 Princeton Univ. Press 
 Princeton (1978) 
 
\end{thebibliography}
\end{document}